\documentclass[12pt]{amsart}
\usepackage{amsmath}
\usepackage{amssymb}
\usepackage{amscd}
\usepackage[mathscr]{eucal}
\usepackage{latexsym}

\newtheorem{thm}{Theorem}[section]
\newtheorem{lem}[thm]{Lemma}

\newtheorem{prop}[thm]{Proposition}

\theoremstyle{definition}

\newtheorem{rem}[thm]{Remark}

\newtheorem{defn}[thm]{Definition}

\newtheorem{ex}[thm]{Example}

\theoremstyle{remark}

\numberwithin{equation}{section}

\def\R{\text{\rm R}}

\def\Gal{\text{\rm Gal}}
\def\Hom{\text{\rm Hom}}

\def\Ker{\text{\rm Ker}\,}

\def\lim{\text{\rm lim}}

\begin{document}

\title[Galois cohomology]
{Galois cohomology of a $p$-adic field via $(\Phi,\Gamma)$-modules in the imperfect residue field case}

\author[Kazuma Morita]{Kazuma Morita}
\address{Department of Mathematics, Faculty of Science,
Kyoto University, Kyoto 606-8502, Japan}
\email{morita@math.kyoto-u.ac.jp}

\date{\today}

\maketitle
\begin{abstract}
For a $p$-adic local field $K$ with perfect residue field, L. Herr constructed a complex which computes the Galois cohomology
of a $p$-torsion representation of the absolute Galois group of $K$ by using the theory of $(\Phi,\Gamma)$-modules. 
We shall generalize his work to the imperfect residue field (the residue field has a finite $p$-basis) case.
\end{abstract}
\section{Introduction}
In this article, $K$ denotes a complete discrete valuation field of characteristic $0$ with residue field $k$ of characteristic $p>0$ such that $[k:k^p]=p^{n}<\infty$.
Assume that $K$ contains a primitive $p$-th root of unity if $p\not=2$ and a primitive $4$-th root of unity if $p=2$.
Choose an algebraic closure $\overline{K}$ of $K$ and put
$G_{K}=\Gal(\overline{K}/K).$
By a $p$-torsion $G_{K}$-representation, we mean a $\mathbb{Z}_{p}$-module of finite length endowed with a continuous action of $G_{K}$.
Let $\text{\bf{Rep}}_{p-\text{tor}}(G_{K})$ denote the category of $p$-torsion $G_{K}$-representations. Let $V$ be a $p$-torsion $G_{K}$-representation.
In the case $n=0$ (i.e. $k$ is a perfect field), Herr [H1] obtained a presentation of the Galois cohomology $H^{*}(G_{K},V)$ in terms of the $(\Phi,\Gamma_{K})$-module
$D(V)$ associated to $V$ in the sense of Fontaine [F].

Now, let $n$ be arbitrary. The purpose of this paper is to give a presentation of $H^{*}(G_{K},V)$ in terms of the $(\Phi,\Gamma_{K})$-module (defined in this paper) associated
to $V$ (Theorem $1.1$). Our $\Gamma_{K}$ is non-commutative if $n\geq 1$. 

Fix a lifting $(b_{i})_{1\leq i\leq n}$ of a $p$-basis of $k$ in $\mathscr{O}_{K}$ (the ring of integers of $K$), and for each $m\geq 1$ and $1\leq i\leq n$, fix a $p^{m}$-th root 
$b_{i}^{1/p^{m}}$ of $b_{i}$ in $\overline{K}$ satisfying $(b_{i}^{1/p^{m+1}})^{p}=b_{i}^{1/p^{m}}$.
Put $K^{(\prime)}=\cup_{m\geq 0}K(b_{i}^{1/p^{m}}, 1\leq i\leq n)$ and $K_{\infty}^{(\prime)}=\cup_{m\geq 0}K^{(\prime)}(\zeta_{p^{m}})$ where 
$\zeta_{p^m}$ denotes a primitive $p^{m}$-th of unity in $\overline{K}$ such that $\zeta_{p^{m+1}}^{p}=\zeta_{p^{m}}$.
The field $K^{(\prime)}$ depends on the choice of a lifting of a $p$-basis of $k$ in $\mathscr{O}_{K}$, but the field $K_{\infty}^{(\prime)}$ doesn't.
Let $K'$ denote the $p$-adic completion of $K^{(\prime)}$. Choose an algebraic closure $\overline{K'} \ (\supset \overline{K})$ of $K'$. Put $K_{\infty}^{\prime}=\cup_{m\geq 0}K^{\prime}(\zeta_{p^{m}})$ in $\overline{K^{\prime}}$.
These fields $K'$ and $K'_{\infty}$ depend on the choice of a lifting of a $p$-basis of $k$ in $\mathscr{O}_{K}$. Put
$\Gamma_{K}=\Gal(K_{\infty}^{(\prime)}/K)$ and $\Gamma_{K^{\prime}}=\Gal(K_{\infty}^{\prime}/K^{\prime}).$
Then, $\Gamma_{K^{\prime}}$ is isomorphic to an open subgroup of $\mathbb{Z}_{p}^{*}$ via the cyclotomic character 
$\chi:\Gamma_{K^{\prime}}\rightarrow \mathbb{Z}_{p}^{*}$ and $\Gamma_{K}$ is isomorphic to the semi-direct product $\Gamma_{K^{\prime}}\ltimes \mathbb{Z}_{p}^{\oplus n}$
where $\Gamma_{K^{\prime}}$ acts on $\mathbb{Z}_{p}^{\oplus n}$ via $\chi$ (see Section $3$). The group $\Gamma_{K}$ is non-commutative if $n\geq 1$.
Since $K^{\prime}$ has perfect residue field which we denote $k'=k^{p^{-\infty}}$, we can apply the theory of Fontaine [F] to obtain the $(\Phi,\Gamma_{K^{\prime}})$-module $D(V)$ for a $p$-torsion $G_{K}$-representation $V$. 
Then, $D(V)$ is equipped with a Frobenius operator $\phi:D(V)\rightarrow D(V)$ and also with a continuous action of $\Gamma_{K}$ (not only $\Gamma_{K^{\prime}}$) which commutes with $\phi$.
With these actions , $D(V)$ becomes an object of the category $ \text{\bf{ $\Phi\Gamma$M}}_{\mathbb{A}_{K},\Gamma_{K}}^{\text{\'et,$p$-tor}}$ of torsion \'etale $(\Phi,\Gamma_{K})$-modules which we will define
(see Section $2$) by imitating the definition of the category of torsion \'etale $(\Phi,\Gamma_{K^{\prime}})$-modules by Fontaine ([F], p273, 3.3.2).
Then, we shall obtain an equivalence of categories between 
$$\text{\bf{Rep}}_{\text{$p$-tor}}(G_{K})\quad\text{and}\quad \text{\bf{$\Phi\Gamma$M}}_{\mathbb{A}_{K},\Gamma_{K}}^{\text{\'et,$p$-tor}}$$ 
which is a generalization of the equivalence of Fontaine ([F], p274, 3.4.3) to the imperfect residue field case (for details, see Theorem $2.7$ in Section $2$).
By using this $D(V)$, we will construct a complex $C_{\phi,\Gamma_{K}}(D(V))$ in Section $3$. Our main result is the following.
\begin{thm}
With notations as above, the group $H^{i}(G_{K},V)$ is canonically isomorphic to the $i$-th cohomology group of the complex $C_{\phi,\Gamma_{K}}(D(V))$ for all $i$.
This isomorphism is functorial in $V$.
\end{thm}
Our proof of the main theorem is a little different from the method of Herr. In the case $n=0$, he considered an ``effaceable'' property of the complex $C_{\phi,\Gamma_{K}}(D(V))$,
whereas our method is to construct a free resolution of the $\mathbb{Z}_{p}[[\Gamma_{K}]]$-module $\mathbb{Z}_{p}$. 

This paper is organized as follows. In Section $2$, we shall review the theory of $(\Phi,\Gamma)$-modules, which is due to J.-M. Fontaine [F] in the perfect residue field case.
We shall construct a theory of  $(\Phi,\Gamma)$-modules in the imperfect residue field case (F. Andreatta [A] constructs a more general and finer theory of $(\Phi,\Gamma)$-modules).
In Section $3$, for $M\in \text{\bf{$\Phi\Gamma$M}}_{\mathbb{A}_{K},\Gamma_{K}}^{\text{\'et,$p$-tor}}$, we shall construct the complexes $C_{\Gamma_{K}}(M)$ and $C_{\phi,\Gamma_{K}}(M)$ which are to be used in the main
theorem. In Section $4$, we shall construct a free resolution of $\mathbb{Z}_{p}$ in the category of left $\mathbb{Z}_{p}[[\Gamma_{K}]]$-modules. In Section $5$, we shall prove that the cohomology group of $C_{\phi,\Gamma_{K}}(D(V))$ coincides with the Galois cohomology
$H^{*}(G_{K},V)$.
\section{The theory of $(\Phi,\Gamma)$-modules}
Let $k'$ denote the perfect residue field of $K^{\prime}$ as in Section $1$. Put $F(k')=W(k')[p^{-1}]$ where $W(k')$ denotes the ring of Witt vectors with coefficients in $k'$.
Now, we apply the theory of $(\Phi,\Gamma)$-modules of Fontaine to $K^{\prime}$.
Since $K^{(\prime)}$ is a Henselian discrete valuation field , we have an isomorphism $G_{K'}=\Gal(\overline{K'}/K')\simeq G_{K^{(\prime)}}=\Gal(\overline{K}/K^{(\prime)}) \ (\subset G_{K})$. With this isomorphism, we identify $G_{K^{\prime}}$ with a subgroup of $G_{K}$.
We have a bijective map from the set of finite extensions of $K^{(\prime)}$ contained in $\overline{K}$ to the set of finite extensions of $K'$ contained in $\overline{K'}$ defined by 
$L\rightarrow LK'$. Furthermore, $LK'$ is the $p$-adic completion of $L$. Hence, we have an
isomorphism of rings
$$\mathscr{O}_{\overline{K}}/p^{n}\mathscr{O}_{\overline{K}}\simeq \mathscr{O}_{\overline{K'}}/p^{n}\mathscr{O}_{\overline{K'}}$$
where $\mathscr{O}_{\overline{K}}$ and $\mathscr{O}_{\overline{K'}}$ denote the rings of integers of $\overline{K}$ and $\overline{K'}$. Thus, the $p$-adic
completion of $\overline{K}$ is isomorphic to the $p$-adic completion of $\overline{K'}$, which we will write
$\mathbb{C}_{p}$. Put
$$\widetilde{\mathbb{E}}=\varprojlim
\mathbb{C}_{p}=\verb+{+(x^{(0)},x^{(1)},\cdots)\verb+|+(x^{(i+1)})^{p}=x^{(i)},x^{(i)}\in\mathbb{C}_{p}\verb+}+$$
and $\widetilde{\mathbb{E}}^{+}$ denotes the set of $x=(x^{(i)})\in \widetilde{\mathbb{E}}$ such that $x^{(0)}\in\mathscr{O}_{\mathbb{C}_{p}}$ (where $\mathscr{O}_{\mathbb{C}_{p}}$ denotes
the ring of integers of $\mathbb{C}_{p}$). For two elements $x = (x^{(i)})$ and $y = (y^{(i)})$ of $\widetilde{\mathbb{E}}$, define
their sum and product by $(x + y)^{(i)} = \lim_{j\rightarrow \infty}(x^{(i+j)} + y^{(i+j)})^{p^j}$ and $(xy)^{(i)}=
x^{(i)}y^{(i)}$. Let $\epsilon = (\epsilon^{(i)})$ denote an element of $\widetilde{\mathbb{E}}$ such that $\epsilon^{(0)} = 1$ and $\epsilon^{(1)}\not=1$.
Then, $\widetilde{\mathbb{E}}$ is a field of characteristic $p > 0$ ($\widetilde{\mathbb{E}}^{+}$ is a subring of $\widetilde{\mathbb{E}}$) and is the completion of an algebraic closure
of $k'((\epsilon-1))$ for the valuation defined by $v_{\mathbb{E}}(x) = v_{p}(x^{(0)})$ where $v_p$ denotes the
$p$-adic valuation of $\mathbb{C}_{p}$ normalized by $v_{p}(p) = 1$.
\begin{ex}  With this valuation, we have
$$v_{\mathbb{E}}(\epsilon - 1) = \lim_{n\rightarrow \infty}v_{p}((\epsilon^{(n)} - 1)^{p^n})=\frac{p}{p-1}.$$
\end{ex}
The field $\widetilde{\mathbb{E}}$ is equipped with an action of a Frobenius $\sigma$ and a continuous action
of the Galois group $G_K$ with respect to the topology defined by the valuation $v_{\mathbb{E}}$.
Put $\mathbb{E}_{F(k')}=k'((\epsilon-1))$ and define $\mathbb{E}$ to be the separable closure of $\mathbb{E}_{F(k')}$ in
$\widetilde{\mathbb{E}}$. Define $H_K =\Gal(\overline{K'}/K_{\infty}')$ which is isomorphic to the subgroup $G_{K_{\infty}^{(\prime)}}=\Gal(\overline{K}/K_{\infty}^{(\prime)})$ of $G_K$. From now on, we identify $H_K$ with $G_{K_{\infty}^{(\prime)}}$. If we put $\mathbb{E}_{K}=
\mathbb{E}^{H_K}$ and define $G_{\mathbb{E}_K}$ to be the Galois group of $\mathbb{E}/\mathbb{E}_K$, the action of $G_{K'}$ on
$\mathbb{E}$ induces the canonical isomorphism $H_K \simeq  G_{\mathbb{E}_K}$ by the theory of the field of
norms ([FW], [W]). Let $\pi$ denote $[\epsilon]-1$. Put $\widetilde{\mathbb{A}}=W(\widetilde{\mathbb{E}})$ ($\widetilde{\mathbb{E}}$ is a perfect field)
and $\widetilde{\mathbb{A}}^{+}=W(\widetilde{\mathbb{E}}^{+})$. The ring $\widetilde{\mathbb{A}}$ is endowed with the topology whose fundamental system of the
neighborhoods of 0 has the form $\pi^{k}\widetilde{\mathbb{A}}^{+}+p^{n+1}\widetilde{\mathbb{A}}$ for $k,n\in\mathbb{N}$. This topology coincides
with the product topology defined by the application $\widetilde{\mathbb{A}}\rightarrow  \widetilde{\mathbb{E}}^{\mathbb{N}}: x \mapsto (x_{k})_{k\in\mathbb{N}}$.
The continuous action of $G_K$ on $\widetilde{\mathbb{E}}$ induces the continuous action of $G_K$ on $\widetilde{\mathbb{A}}$
which commutes with the Frobenius $\sigma$. Let $\mathbb{A}_{F(k')}$ be the $p$-adic completion of
$W(k')[[\pi]][\pi^{-1}]$ contained in $\widetilde{\mathbb{A}}$. This ring is a complete discrete valuation ring
with the residue field $\mathbb{E}_{F(k')}$. Let $\mathbb{A}$ be the $p$-adic completion of the maximal
unramified extension of $\mathbb{A}_{F(k')}$ in $\widetilde{\mathbb{A}}$ which has the residue field $\mathbb{E}$. The ring $\mathbb{A}$ is
equipped with an action of the Galois group $G_K$ and of the Frobenius $\sigma$ induced
from those of $\widetilde{\mathbb{E}}$. Put $\mathbb{A}_{K} = \mathbb{A}^{H_{K}}$.

For all $V \in \text{\bf{Rep}}_{p-\text{tor}}(G_{K'})$, we can associate the $(\Phi,\Gamma_{K'})$-module over $\mathbb{A}_{K}$
$$D(V) = (\mathbb{A}\otimes_{\mathbb{Z}_{p}}V)^{H_{K}}.$$
It is equipped with the residual action of $\Gamma_{K'}\simeq G_{K'}/H_{K}$ and the Frobenius $\phi_{D(V)}$
induced by that on $\mathbb{A}$. The module $D(V)$ is a torsion \'etale $(\Phi,\Gamma_{K'})$-module over
$\mathbb{A}_{K}$ ([F], p274, 3.4.2).

Conversely, to a torsion \'etale $(\Phi,\Gamma_{K'})$-module $M$ over $\mathbb{A}_{K}$, we can associate a
$p$-torsion representation of $G_{K'}$ as follows
$$(*)\quad V(M)=(\mathbb{A}\otimes_{\mathbb{A}_{K}}M)^{\sigma\otimes\phi_{M}=1}\in\text{\bf{Rep}}_{p-\text{tor}}(G_{K'}).$$
Let $\text{\bf{$\Phi\Gamma$M}}_{\mathbb{A}_{K},\Gamma_{K'}}^{\text{\'et,$p$-tor}}$
denote the category of torsion \'etale $(\Phi,\Gamma_{K'})$-modules in the sense
of Fontaine ([F], p273, 3.3.2). By the two constructions above, Fontaine proved
the following ([F], p274, 3.4.3).
\begin{thm} The functor $D$ gives an equivalence between the two categories
$$\text{\bf{Rep}}_{p-\text{tor}}(G_{K'})\quad \text{and}\quad \text{\bf{$\Phi\Gamma$M}}_{\mathbb{A}_{K},\Gamma_{K'}}^{\text{\'et,$p$-tor}}$$
The functor $V$ is a quasi-inverse of $D$.
\end{thm}
Define $(\Phi,\Gamma_{K})$-modules as follows.
\begin{defn} A torsion $(\Phi,\Gamma_{K})$-module over $\mathbb{A}_{K}$ is an $\mathbb{A}_{K}$-module $M$ of finite
length equipped with
 \begin{enumerate}
\item a $\sigma$-semi-linear map (which we call a Frobenius operator)
$$\phi=\phi_{M}:M\rightarrow M$$
 \item a continuous semi-linear action of $\Gamma_{K}$ which commutes with $\phi$.
\end{enumerate}
In addition, we call $M$ an \'etale $(\Phi,\Gamma_{K})$-module if it is generated by the image of
$\phi$ as an $\mathbb{A}_{K}$-module. Let $\text{\bf{$\Phi\Gamma$M}}_{\mathbb{A}_{K},\Gamma_{K}}^{\text{\'et,$p$-tor}}$
denote the category which consists of

$\circ $ objects: torsion \'etale $(\Phi,\Gamma_{K})$-modules over $\mathbb{A}_{K}$

$\circ $ morphisms: $\mathbb{A}_{K}$-linear morphisms which commute with $\phi$ and the action of
$\Gamma_{K}$.
\end{defn}
\begin{rem} Put $\mathbb{A}^{+}_{K}=\mathbb{A}_{K}\cap \widetilde{\mathbb{A}}^{+}$. If we fix a lifting $T_{K}$ of the prime element
of $\mathbb{E}_{K}$ in $\mathbb{A}_{K}^{+}$, we have $\mathbb{A}_{K}^{+}=W(k')[[T_{K}]]$. Let $M$ be a finitely generated $\mathbb{A}_{K}/p^{n}$-module.
Fix a finitely generated sub-$\mathbb{A}_{K}^{+}/p^{n}$-module $M^0$ of $M$ such that $M$ is
generated by $M^0$ over $\mathbb{A}_{K}/p^{n}$. The module $M$ is endowed with the topology 
such that the family of submodules $\verb+{+T_{K}^{m}M^{0}\verb+}+_{m\geq 1}$ is a fundamental system of neighborhoods of $0$. 
This topology is
independent of the choice of $M^0$. Furthermore, since $\mathbb{A}_{K}/p^n$ is Noetherian and
complete for the $T_K$-adic topology, $T_{K}^{-N}M^{0}$ is complete for the $T_K$-adic topology.
We may use the family of submodules $\verb+{+\pi^{m}M^{0}\verb+}+_{m\geq 1}$ instead of  $\verb+{+T_{K}^{m}M^{0}\verb+}+_{m\geq 1}$ to define the same topology.
\end{rem}
Consider a $p$-torsion $G_K$-representation $V \in \text{\bf{Rep}}_{p-\text{tor}}(G_K)$ and $D(V)$. Since
the Galois group $G_K$ acts on $\mathbb{A}\otimes_{\mathbb{Z}_{p}}V$ and we have $D(V)=(\mathbb{A}\otimes_{\mathbb{Z}_{p}}V)^{H_K}$, the
quotient $\Gamma_{K}\simeq G_{K}/H_{K}$ and $\phi$ act on $D(V)$ commuting with each other. This
means that $D(V)$ becomes an object of $\text{\bf{$\Phi\Gamma$M}}_{\mathbb{A}_{K},\Gamma_{K}}^{\text{\'et,$p$-tor}}$
. The continuous action of
$G_K$ on $\mathbb{A}\otimes_{\mathbb{Z}_{p}}V$ induces the continuous action of $\Gamma_{K}$ on $D(V)$ as follows. Let
$L$ be a finite Galois extension of $K$ contained in $\overline{K}$ such that the action of
$G_{L} = \Gal(\overline{K}/L)$ on $V$ is trivial. Fix $n \in \mathbb{N}$ such that $p^{n}V = 0$. Then, we
have $D(V)=(\mathbb{A}_{L}/p^{n}\otimes_{\mathbb{Z}_{p}}V)^{H_{K}}$. The following two topologies of $D(V)$ coincide
\begin{enumerate}
\item the topology defined in Remark 2.4
\item the induced topology as a subspace of $\mathbb{A}_{L}/p^{n}\otimes_{\mathbb{Z}_{p}}V$ whose topology is
defined in Remark 2.4.
\end{enumerate}
(Proof: There exists an $\mathbb{A}_{L}/p^{n}$-linear isomorphism
$$\mathbb{A}_{L}/p^{n}\otimes_{\mathbb{A}_{K}/p^{n}}D(V)\simeq \mathbb{A}_{L}/p^{n}\otimes_{\mathbb{Z}_{p}}V.$$
Fix a finitely generated sub-$\mathbb{A}_{K}^{+}/p^{n}$-module $M^{0}$ of $D(V)$ such that $D(V)$ is generated
by $M^0$ over $\mathbb{A}_{K}/p^n$. Let $M^{0}_{L}$
be the sub-$\mathbb{A}_{L}^{+}/p^{n}$-module of $\mathbb{A}_{L}/p^{n}\otimes_{\mathbb{Z}_{p}}V$
generated by $M^0$. Since the morphism $\mathbb{A}^{+}_{
K}/p^{n}\rightarrow  \mathbb{A}^{+}_{L}/p^n$ is finite flat, the morphism
$\mathbb{A}_{L}^{+}/p^{n}\otimes_{\mathbb{A}_{K}^{+}/p^{n}}M^{0}\rightarrow M_{L}^{0}$
is an isomorphism. Thus, the inverse image of $\pi^{m}M_{L}^{0}$
by the map $D(V)\rightarrow  \mathbb{A}_L/p^{n}\otimes_{\mathbb{Z}_{p}}V$ is $\pi^{m}M^{0}$.)
\begin{rem} Let $L$ be a finite Galois extension of $K$ contained in $\overline{K}$. Let $M$
be a finitely generated $\mathbb{A}_{L}/p^n$-module endowed with a continuous and semi-linear
action of $\Gal(K_{\infty}^{(\prime)}L/K)$. Fix $M^{0}$ as in Remark 2.4. Let $M^1$ be the sub-$\mathbb{A}^{+}
_L/p^n$-module
of $M$ generated by $g(M^0)$ $(g \in \Gal(K_{\infty}^{(\prime)}L/K))$. Since $\Gal(K_{\infty}^{(\prime)}L/K)$
is compact, $M^1$ is also a finitely generated sub-$\mathbb{A}^{+}_
L/p^n$-module. By construction, $M^{1}$ is stable under the action of $\Gal(K_{\infty}^{(\prime)}L/K)$. Then, for $N,m \in \mathbb{N}$,
$\pi^{-N}M^{1}/\pi^{m}M^{1}$ becomes a discrete $\mathbb{Z}_{p}[[\Gal(K_{\infty}^{(\prime)}L/K)]]$-module. 
Since $\pi^{-N}M^{1}$ is
complete for the $\pi$-adic topology, $\pi^{-N}M^{1}$ has the structure of $\mathbb{Z}_{p}[[\Gal(K_{\infty}^{(\prime)}L/K)]]$-
module. Thus, $M$ is equipped with the structure of $\mathbb{Z}_{p}[[\Gal(K_{\infty}^{(\prime)}L/K)]]$-module.
For another sub-$\mathbb{A}^{+}_{L}/p^{n}$-module $M^2$ of $M$ stable under the action of $\Gal(K_{\infty}^{(\prime)}L/K)$ such that $M$ is generated by $M^2$ over $\mathbb{A}_{L}/p^n$, we
can find integers $N_1 \leq N_2$ such that $\pi^{N_{1}}M^{1}\subset M^{2}\subset \pi^{N_{2}}M^{1}$, therefore, the
structure of $\mathbb{Z}_{p}[[\Gal(K_{\infty}^{(\prime)}L/K)]]$-module on $M$ is independent of the choice of
$M^1$. With this, the action of $\Gamma_{K}$ on $D(V)$ naturally extends to the action of $\mathbb{Z}_{p}[[\Gamma_{K}]]$.
\end{rem}
\begin{rem}Let $L$ be a finite Galois extension of $K$ contained in $\overline{K}$ such
that the action of $G_L$ on $V$ is trivial. Fix $n \in \mathbb{N}$ such that $p^{n}V = 0$. For a
finite Galois extension of $K$ such that $L \subset L'\subset \overline{ K}$, $\mathbb{A}_{L'}\otimes_{\mathbb{Z}_{p}}V$ is a finitely
generated $\mathbb{A}_{L'}/p^n$-module endowed with a continuous and semi-linear action of
$\Gal(K_{\infty}^{(\prime)}L'/K)$. Remark 2.5 says that the action of $\Gal(K_{\infty}^{(\prime)}L'/K)$ on $\mathbb{A}_{L'}\otimes_{\mathbb{Z}_{p}}V$
naturally extends to the action of $\mathbb{Z}_p[[\Gal(K_{\infty}^{(\prime)}L'/K)]]$. Thus, the action of $G_K$
on $\mathbb{A}\otimes_{\mathbb{Z}_{p}}V= \varinjlim _{L'}(\mathbb{A}_{L'}\otimes_{\mathbb{Z}_{p}}V)$ naturally extends to the action of $\mathbb{Z}_{p}[[G_{K}]]$.
Then, the canonical injection $D(V)\rightarrow \mathbb{A}\otimes_{\mathbb{Z}_{p}}V$ is compatible with the action
of $\mathbb{Z}_{p}[[G_{K}]]$. (Proof: Let $L$, $M^0$, $M^{0}_{L}$
be as in the proof of the coincidence of the two
topologies of $D(V )$ before Remark 2.5. We can assume that $M^0$ is endowed
with a continuous and semi-linear action of $\Gamma_K$ (see Remark 2.5). Since the
morphism $\mathbb{A}^{+}_{K}/p^{n}\rightarrow \mathbb{A}_{L}^{+}/p^{n}$ is finite flat, we have a morphism of discrete $\mathbb{Z}_{p}[[G_K]]$-
modules
$\pi^{-N}M^{0}/\pi^{m}M^{0}\rightarrow (\mathbb{A}_{L}^{+}/p^{n})\otimes_{\mathbb{A}_{K}^{+}/p^{n}}(\pi^{-N}M^{0}/\pi^{m}M^{0})\simeq\pi^{-N}M_{L}^{0}/\pi^{m}M_{L}^{0}$. 
By taking the inverse limit for $m$, we obtain the morphism of $\mathbb{Z}_{p}[[G_K]]$-modules
$\pi^{-N}M^{0}\rightarrow \pi^{-N}M_{L}^{0}$.)
\end{rem}
Conversely, to a torsion \'etale ($\Phi,\Gamma_{K}$)-module $M$ over $\mathbb{A}_K$, we can associate a
$p$-torsion representation of $G_K$ as follows (see ($*$))
$$V (M) = (\mathbb{A}\otimes_{\mathbb{A}_{K}}M)^{\sigma\otimes\phi_{M}=1}\in.$$
The continuous action of $G_K$ on $\mathbb{A}\otimes_{\mathbb{A}_{K}}M$ induces the continuous action of $G_K$
on $V (M).$ Here, we give $V (M)$ the induced topology as a subspace of $\mathbb{A}\otimes_{\mathbb{A}_{K}}M$.
Since the topology of $\mathbb{A}\otimes_{\mathbb{A}_{K}}M$ is Hausdorff and $V (M)$ is finite, the induced
topology on $V (M)$ is discrete.

An imperfect residue field version of Fontaine's theorem is the following (cf.
Theorem 2.2).
\begin{thm} The functor $D$ gives an equivalence between the two categories
$$\text{\bf{Rep}}_{\text{$p$-tor}}(G_{K})\quad\text{and}\quad \text{\bf{$\Phi\Gamma$M}}_{\mathbb{A}_{K},\Gamma_{K}}^{\text{\'et,$p$-tor}}$$ 
The functor $V$ is a quasi-inverse of $D.$
\end{thm}
\begin{proof} For $M\in \text{\bf{$\Phi\Gamma$M}}_{\mathbb{A}_{K},\Gamma_{K}}^{\text{\'et,$p$-tor}}$, the natural morphism
$$\mathbb{A}\otimes_{\mathbb{Z}_{p}}V(M)\rightarrow \mathbb{A}\otimes_{\mathbb{A}_{K}}M$$
induces a morphism $D(V (M)) \rightarrow M$ and this morphism is an isomorphism ([F],
p258, 1.2.6). Conversely, for $N \in \text{\bf{Rep}}_{\text{$p$-tor}}(G_{K})$, the natural morphism
$$\mathbb{A}\otimes_{\mathbb{A}_{K}}D(N)\rightarrow \mathbb{A}\otimes_{\mathbb{Z}_{p}}N$$ 
induces a morphism $V (D(N)) \rightarrow N$ and this morphism is an isomorphism ([F],
p258, 1.2.4). 
\end{proof}

\section{Main theorem}
We will give a presentation of $H^{*}(G_{K},V)$ in terms of $D(V)$. Recall that we fixed a $p^{m}$-th root $b_{i}^{1/p^{m}}$ of $b_{i}$ in Introduction. Fix a $p^{m}$-th root $\zeta_{p^{m}}$ of unity such that $\zeta_{p^{m+1}}^{p}=\zeta_{p^{m}}$.  
Fix a topological generator $\gamma$ of $\Gamma_{K'}\subset \Gamma_{K}$ and define $\beta_{i}\in \Gamma_{K}$ ($1\leq i\leq n$) by
$$\beta_{i}(b_{i}^{1/p^m})=b_{i}^{1/p^m}\zeta_{p^m},\mspace{5mu} \beta_{i}(b_{j}^{1/p^m})=b_{j}^{1/p^m}\mspace{5mu}(j\not=i)\mspace{5mu}\text{and}\mspace{5mu}\beta_{i}(\zeta_{p^m})=\zeta_{p^m}.$$
Define $l\in\mathbb{Z}_{p}^{*}$ by 
$$\gamma(\zeta_{p^m})=\zeta_{p^m}^{l}.$$
These topological generators $(\gamma,\beta_{1},\ldots,\beta_{n})$ define the isomorphism $\Gamma_{K}\simeq \Gamma_{K^{\prime}}\ltimes \mathbb{Z}_{p}^{\oplus n}$ ($\beta_{i}\mapsto$ the topological generator of $i$-th component of $\mathbb{Z}_{p}$).     
Let $\Lambda$ denote $\mathbb{Z}_{p}[[\Gamma_{K}]]$ in what follows.
Define elements of $\Lambda$ as follows
$$\omega_{i}=\beta_{i}-1\mspace{15mu}\text{and}\mspace{15mu}\tau_{S}=(\Pi_{i\in S}\frac{\beta_{i}-1}{\beta_{i}^{l}-1})\gamma-1.$$
Recall that $D(V )$ is naturally equipped with the action of $\Lambda$ (Remark 2.5).
Since $(\beta^{l}_{i}-1)(\beta_{i}-1)^{-1}=\verb+{+(1+\omega_{i})^{l}-1\verb+}+\omega_{i}^{-1}\in l+\omega_{i}\mathbb{Z}_{p}[[\omega_{i}]]$ and $l\in \mathbb{Z}_{p}^{*}$, we have $(\beta_{i}^{l}-1)(\beta_{i}-1)^{-1}\in\mathbb{Z}_{p}[[\omega_{i}]]^{*}$.
\begin{enumerate}
\item 
{\bf The complex $C_{\Gamma_K}(D(V))$}

To a $p$-torsion representation $V$ of $G_{K}$, define the complex $C_{\Gamma_K}$ $(D(V))$ to be 
$$0\longrightarrow D(V)^{X(0)}\mathop{\longrightarrow}^{{d}^0} D(V)^{X(1)}\mathop{\longrightarrow}^{{d}^1}\cdots \\
\mathop{\longrightarrow}^{{d}^{i-1}} D(V)^{X(i)}\mathop{\longrightarrow}^{{d}^{i}}\cdots \\
\mathop{\longrightarrow}^{{d}^{n}} D(V)^{X(n+1)}\longrightarrow 0.
$$
(The proof of $d^{i}\circ d^{i-1}=0$ follows from the presentation of $C_{\Gamma_{K}}(D(V))$ in terms of $C_{\Lambda}$ in Section $5$.)

Here
\begin{enumerate}
\item For a finite set $X$, we define $D(V)^{X}=\oplus_{S\in X}D(V)$.
\item $X(i)$ denotes the set of all subsets of $\verb+{+ 0,\cdots,n\verb+}+$ of order $i$. Notice that the order of $X(i)$ is ${n+1\choose i}$.
We define the degree of $D(V)^{X(0)}$ to be $0$.
\item For $S\in X(i)$ and $T\in X(i+1)$, the $(S,T)$-component $d^{i}(S,T)$ of 
$d^{i}:D(V)^{X(i)}\rightarrow D(V)^{X(i+1)}$
is defined as follows.

\vspace{10pt}

(A) If $S\not\subset T$, $d^{i}(S,T)=0$.

\vspace{10pt}

(B) If $S\subset T$, put $\verb+{+j\verb+}+=T\verb+\+S$.

\vspace{10pt}

    $\circ$ If $j=0$, $d^{i}(S,T)=\tau_{S}$.

\vspace{10pt}

    $\circ$ If $j\not=0$, $d^{i}(S,T)=(-1)^{a(S,j)}\omega_{j}$ 
where $a(S,j)=\sharp\verb+{+x\in S;x\leq j\verb+}+$.
\end{enumerate}

\vspace{10pt}

\item
{\bf The complex $C_{\phi,\Gamma_K}(D(V))$}

Define the complex $C_{\phi,\Gamma_K}(D(V))$ by
$$C_{\phi,\Gamma_{K}}(D(V))=\mspace{5mu}\text{the mapping fiber of}\ C_{\Gamma_{K}}(D(V))\mathop{\rightarrow}^{\rho} C_{\Gamma_{K}}(D(V))$$
where $\rho=\phi-1$.
The complex $C_{\phi,\Gamma_{K}}(D(V))$ has the following form
\begin{align*}
0\longrightarrow D(V)^{\oplus {n+2\choose 0}} \mathop{\longrightarrow}^{d^0} D(V)^{\oplus {n+2\choose 1}}&  \mathop{\longrightarrow}^{d^1} \cdots\\ 
\mathop{\longrightarrow}^{d^{i-1}}& D(V)^{\oplus {n+2\choose i}}  \mathop{\longrightarrow}^{d^{i}} \cdots 
\mathop{\longrightarrow}^{d^{n+1}} D(V)^{\oplus {n+2\choose n+2}} \longrightarrow 0.
\end{align*}
(Here, define the degree of $D(V)^{\oplus {n+2\choose 0}}$ to be $0$.)
\end{enumerate}

Our main result is the following.

\begin{thm}With notations as above, the group $H^i(G_K,V)$ is canonically isomorphic to the $i$-th cohomological group of the complex $C_{\phi,\Gamma_{K}}(D(V))$ for all $i$.
This isomorphism is functorial in $V$.
\end{thm}
It follows that the cohomological dimension of $K$ is $n+2$.
\begin{ex}
\begin{enumerate}
\item The case $n=0$ (i.e. the residue field $k$ is perfect)

In this case, the complex $C_{\phi,\Gamma_K}(D(V))$ is given by
$$0\longrightarrow D(V)\mathop{\longrightarrow}^{d^0}D(V)\oplus D(V)\mathop{\longrightarrow}^{d^1}D(V)\longrightarrow0$$

$\circ$ $d^0(x)=(\rho(x),\tau(x))$,

\vspace{10pt}

$\circ$ $d^1(x,y)=(\rho(y)-\tau(x))$.

\vspace{10pt}

Here, $\rho=\phi-1$ and $\tau=\tau_{\varnothing}=\gamma-1$.
This is the complex constructed by Herr.
\item The case $n=1$ (i.e. the residue field $k$ is imperfect and $[k:k^{p}]=p$)

In contrast to the example $(1)$, there is an action of $\omega_1=\beta_1-1$.
Therefore, we have a more complicated complex than before.
\begin{align*}0\longrightarrow D(V)\mathop{\longrightarrow}^{d^0}D(V)&\oplus D(V)\oplus D(V)\\
\mathop{\longrightarrow}^{d^1}&D(V)\oplus D(V)\oplus D(V)\mathop{\longrightarrow}^{d^2}D(V)\longrightarrow 0
\end{align*}

$\circ$ $d^0(x)=(\rho(x),\tau(x),\omega_1(x)),$

\vspace{10pt}

$\circ$ $d^1(x,y,z)=(\rho (y)-\tau (x),\rho (z)-\omega_1(x),\tau_{\verb+{+1\verb+}+}(z)-\omega_1(y)),$

\vspace{10pt}

$\circ$ $d^2(x,y,z)=(\rho(z)-\tau_{\verb+{+1\verb+}+}(y)+\omega_1(x)).$

\vspace{10pt}

The appearance of $\tau_{\verb+{+1\verb+}+}$, instead of $\tau$, reflects the non-commutativity of $\Gamma_{K}$.
\end{enumerate}
\end{ex}
\section{Construction of a free resolution of $\mathbb{Z}_p$}
\subsection{Relations in $\Lambda$}
Let  
$$\gamma,\beta_1,\beta_2,\cdots,\beta_{n}$$
be the topological generators of $\Gamma_{K}$ as in the previous section.
We have the following 
relations
\begin{enumerate}
\item $\gamma\beta_i={\beta_i}^{l} \gamma \mspace{5mu}(l\in \mathbb{Z}_{p}^{*})$
\item $\beta_i\beta_j=\beta_j\beta_i$.
\end{enumerate}
For the construction of the complex, consider the following elements of $\Lambda$ (these are introduced in the previous section)
\begin{enumerate}
\item $\tau =\gamma -1$ 
\item $\omega_i=\beta_i-1$ 
\item $W_i={\beta_i}^{l}-1$
\item $\tau_{S}=\underset{i\in S}{\Pi} (\omega_{i} W_{i}^{-1})\gamma-1$ for $S\subset \verb+{+1,\cdots,n\verb+}+$. 

Recall that $\omega_{i}W_{i}^{-1}\in \mathbb{Z}_{p}[[\beta_{i}-1]]\subset \Lambda$.
\end{enumerate}
\begin{rem} Notice that $\tau_{S}=\tau$ if $S=\varnothing$.
\end{rem}
These operators have the following relations:

{\bf Relations (R)}
\begin{enumerate}

\item $\omega_i\omega_j=\omega_j\omega_i$
\item $W_i W_j=W_j W_i$ 
\item $\gamma\omega_{i}=W_{i}\gamma$
\item For $i\in S\subset\verb+{+1,\cdots ,n\verb+}+$, $\tau_{S} \omega_i=\omega_i\tau_{S\verb+\+\verb+{+i\verb+}+}.$
\begin{proof}
$\tau_{S} \omega_i$

$=({\underset{j\in S}{\Pi}}(\omega_{j}W_{j}^{-1})\gamma-1)\omega_i$
$=({\underset{j\in S}{\Pi}}(\omega_{j}W_{j}^{-1})\gamma \omega_i-\omega_i)$

$=({\underset{j\in S}{\Pi}}(\omega_{j}W_{j}^{-1})W_i\gamma-\omega_i)$
$=\omega_i({\underset{j\in S,j\not=i}{\Pi}}(\omega_{j}W_{j}^{-1})\gamma-1)$

$=\omega_i\tau_{S\verb+\+\verb+{+i\verb+}+}.$
\end{proof}
\end{enumerate}
\subsection{Construction of $C_{\Lambda}$}
Consider the following sequence $C_{\Lambda}$ of left $\Lambda$-modules
$$0\longrightarrow \Lambda ^{X(n+1)} \mathop{\longrightarrow}^{d_n} \Lambda^{X(n)} \mathop{\longrightarrow}^{d_{n-1}} \cdots
\mathop{\longrightarrow}^{d_i} \Lambda^{X(i)}\mathop{ \longrightarrow}^{d_{i-1}} \cdots
\mathop{\longrightarrow}^{d_{0}} \Lambda^{X(0)}\longrightarrow 0.$$
Here
\begin{enumerate}
\item For a finite set $X$, define $\Lambda^{X}=\oplus_{S\in X}\Lambda$.
\item $X(i)$ denotes the set of all subsets of $\verb+{+0,1,\cdots,n\verb+}+$ of order $i$.
 Define the degree of $\Lambda^{X(0)}$ to be $0$.
\item For $S\in X(i)$ and $T\in X(i+1)$, the $(S,T)$-component $d_{i}(S,T)$ of $d_{i}:\Lambda^{X(i+1)}\rightarrow \Lambda^{X(i)}$ is defined as follows.

\vspace{10pt}

(A) If $S\not\subset T$, $d_{i}(S,T)(x)=0$.

\vspace{10pt}

(B) If $S\subset T$, put $\verb+{+j\verb+}+=T\verb+\+S$.

\vspace{10pt}

$\circ$ If $j=0$, $d_{i}(S,T)(x)=x\tau_{S}.$

\vspace{10pt}

$\circ$ If $j\not=0$, $d_{i}(S,T)(x)=(-1)^{a(S,j)}x\omega_{j}$ where $a(S,j)=\sharp\verb+{+y\in S;y\leq j\verb+}+$.
\end{enumerate} 
\begin{ex}
In the case of $[k:{k}^p]=p$, we have
$$0\longrightarrow\Lambda\mathop{\longrightarrow}^{d_1}\Lambda^{\oplus 2}\mathop{\longrightarrow}^{d_0}\Lambda\longrightarrow 0$$
Here, $d_0(f,g)=(f\tau+g\omega_1)$ and $d_1(f)=(-f\omega_1, f\tau_{\verb+{+1\verb+}+})$. 
\end{ex}
\begin{lem}
The natural morphism 
$$\underleftarrow{\lim}_{m} \ \mathbb{Z}_{p}[\Gamma_{K'}]/(\Gamma_{K'})^{p^{m}}\otimes_{\mathbb{Z}_{p}}\mathbb{Z}_{p}[[\omega_{1},\ldots,\omega_{n}]]\rightarrow \Lambda=\mathbb{Z}_{p}[[\Gamma_{K}]]$$
is an isomorphism of left $\mathbb{Z}_{p}[[\Gamma_{K'}]]$- and right $\mathbb{Z}_{p}[[\omega_{1},\ldots,\omega_{n}]]$-modules.
\end{lem}
\begin{proof}
For $m\in\mathbb{N}_{>0}$, put $\Gamma_{m}=\Gamma_{K'}/(\Gamma_{K'})^{p^{m}}\ltimes (\mathbb{Z}/p^{m}\mathbb{Z})^{\oplus n}$.
Note that the action of $\Gamma_{K'}$ on $(\mathbb{Z}/p^{m}\mathbb{Z})^{\oplus n}$ factors through the quotient $\Gamma_{K'}/(\Gamma_{K'})^{p^{m}}.$
Then, we have $\mathbb{Z}_{p}[[\Gamma_{K}]]\simeq \underleftarrow{\lim}_{m}\mathbb{Z}_{p}[\Gamma_{m}]$.
The natural homomorphisms of rings $f:\mathbb{Z}_{p}[\Gamma_{K'}/(\Gamma_{K'})^{p^{m}}]\rightarrow \mathbb{Z}_{p}[\Gamma_{m}]$ and $g:\mathbb{Z}_{p}[(\mathbb{Z}/p^{m}\mathbb{Z})^{\oplus n}]\rightarrow \mathbb{Z}_{p}[\Gamma_{m}]$ induce 
the surjection of left $\mathbb{Z}_{p}[\Gamma_{K'}/(\Gamma_{K'})^{p^{m}}]$- and right $\mathbb{Z}_{p}[(\mathbb{Z}/p^{m}\mathbb{Z})^{\oplus n}]$-modules
$$\mathbb{Z}_{p}[\Gamma_{K'}/(\Gamma_{K'})^{p^{m}}]\otimes_{\mathbb{Z}_{p}}\mathbb{Z}_{p}[(\mathbb{Z}/p^{m}\mathbb{Z})^{\oplus n}]\twoheadrightarrow \mathbb{Z}_{p}[\Gamma_{m}]:a\otimes b\mapsto f(a)g(b).$$
Since both sides have the same $\mathbb{Z}_{p}$-rank, it turns out to be an isomorphism.
On the other hand, we have
\footnotesize{\begin{align*}\underleftarrow{\lim}_{m}(\mathbb{Z}_{p}[\Gamma_{K'}/(\Gamma_{K'})^{p^{m}}]\otimes_{\mathbb{Z}_{p}}\mathbb{Z}_{p}[(\mathbb{Z}/p^{m}\mathbb{Z})^{\oplus n}])&\simeq \underleftarrow{\lim}_{m}(\underleftarrow{\lim}_{m'}(\mathbb{Z}_{p}[\Gamma_{K'}/(\Gamma_{K'})^{p^{m}}]\otimes_{\mathbb{Z}_{p}}\mathbb{Z}_{p}[(\mathbb{Z}/p^{m'}\mathbb{Z})^{\oplus n}]))\\                      
                                                                                                                                                                         &\simeq \underleftarrow{\lim}_{m}(\mathbb{Z}_{p}[\Gamma_{K'}/(\Gamma_{K'})^{p^{m}}])\otimes_{\mathbb{Z}_{p}}\mathbb{Z}_{p}[[\omega_{1},\ldots,\omega_{n}]].
\end{align*}}
\normalsize{This completes the proof.}
\end{proof}
\begin{prop}
The sequence 
$$0\longrightarrow\Lambda ^{X(n+1)}\mathop{\longrightarrow}^{d_n}\Lambda ^{X(n)}\mathop{\longrightarrow}^{d_{n-1}}\cdots\mathop{\longrightarrow}^{d_{0}}\Lambda ^{X(0)}\mathop{\longrightarrow}^{Aug} \mathbb{Z}_p\longrightarrow 0$$
gives a free resolution of the left $\Lambda$-module $\mathbb{Z}_p$.
Here, $\mathbb{Z}_{p}$ is equipped with the structure of left $\Lambda$-modules induced from the trivial action of $\Gamma_{K}$.
\end{prop}
\begin{proof}
Consider the following sequence $C_{\omega}$
$$0\longrightarrow \Lambda ^{Y(n)} \mathop{\longrightarrow}^{d'_{n-1}} \Lambda^{Y(n-1)} \mathop{\longrightarrow}^{d'_{n-2}} \cdots
\mathop{\longrightarrow}^{d'_{0}} \Lambda^{Y(0)}\mathop{\longrightarrow}^{Aug} \mathbb{Z}_{p}[[\tau]]\longrightarrow 0.$$
Here
\begin{enumerate}
\item $Y(i)$ denotes the set of all subsets of $\verb+{+1,\cdots,n\verb+}+$ of order $i$.
(Recall that $X(i)$ denotes the set of all subsets of $\verb+{+0,1,\cdots,n\verb+}+$.)
\item For $S\in Y(i)$ and $T\in Y(i+1)$, the $(S,T)$-component $d'_{i}(S,T)$ of $d'_{i}:\Lambda^{Y(i+1)}\rightarrow \Lambda^{Y(i)}$ is defined as follows.

\vspace{10pt}

$\circ$ If $S\not\subset T$, $d'_{i}(S,T)(x)=0$.

\vspace{10pt}

$\circ$ If $S\subset T$, put $\verb+{+j\verb+}+=T\verb+\+S.$

$d'_{i}(S,T)(x)=(-1)^{a(S,j)}x\omega_{j}$ where $a(S,j)=\sharp\verb+{+y\in S;y\leq j\verb+}+$.
\end{enumerate}
Put $\Lambda_{0}=\mathbb{Z}_{p}[[\omega_{1},\ldots,\omega_{n}]].$
Let $K.(\omega_{1},\ldots,\omega_{n})$ be the Koszul complex
$$0\longrightarrow \Lambda ^{Y(n)}_{0} \mathop{\longrightarrow}^{d'_{n-1}} \Lambda^{Y(n-1)}_{0} \mathop{\longrightarrow}^{d'_{n-2}} \cdots
\mathop{\longrightarrow}^{d'_{0}} \Lambda^{Y(0)}_{0}\longrightarrow 0.$$
Since we have the isomorphism $\underleftarrow{\lim}_{m} \ \mathbb{Z}_{p}[\Gamma_{K'}/(\Gamma_{K'})^{p^{m}}]\otimes_{\mathbb{Z}_{p}}\mathbb{Z}_{p}[[\omega_{1},\ldots,\omega_{n}]]\simeq \Lambda=\mathbb{Z}_{p}[[\Gamma_{K}]]$, the sequence 
$$0\longrightarrow \Lambda ^{Y(n)} \mathop{\longrightarrow}^{d'_{n-1}} \Lambda^{Y(n-1)} \mathop{\longrightarrow}^{d'_{n-2}} \cdots
\mathop{\longrightarrow}^{d'_{0}} \Lambda^{Y(0)}\longrightarrow 0$$
is the complex $\underleftarrow{\lim}_{m} \ \mathbb{Z}_{p}[\Gamma_{K'}/(\Gamma_{K'})^{p^{m}}]\otimes_{\mathbb{Z}_{p}}K.(\omega_{1},\ldots,\omega_{n})$, so the sequence $C_{\omega}$ is a resolution of $\mathbb{Z}_{p}[[\Gamma_{K'}]]=\Lambda/(\sum_{i=1}^{n}\Lambda \omega_{i})$ in the category of left $\Lambda$-modules (note that the transition map $\mathbb{Z}_{p}[\Gamma_{K'}/(\Gamma_{K'})^{p^{m}}]\rightarrow \mathbb{Z}_{p}[\Gamma_{K'}/(\Gamma_{K'})^{p^{m'}}]$ ($m\geq m'$) is surjective).
In particular, the sequence $C_{\omega}$ is exact.
Then, consider the following commutative diagram of left $\Lambda$-modules (the commutativity follows from the relation (4)):
$$\begin{CD}
@. @. @. @. @. 0 @. \\
@. @. @. @. @. @VVV @. \\
0 @>>> \Lambda^{Y(n)} @>d'_{n-1}>> \Lambda^{Y(n-1)} @>d'_{n-2}>> \cdots @>d'_{0}>> \Lambda^{Y(0)} @>>> \mathbb{Z}_p[[\tau]] @>>> 0 \\
@. @Vd^{''}_{n}VV @Vd^{''}_{n-1}VV @. @Vd^{''}_{0}VV @V\tau VV @. \\
0 @>>> \Lambda^{Y(n)} @>d'_{n-1}>> \Lambda^{Y(n-1)} @>d'_{n-2}>> \cdots @>d'_{0}>> \Lambda^{Y(0)} @>>> \mathbb{Z}_p[[\tau]] @>>> 0 \\
@. @. @. @. @. @VVV @. \\
@. @. @. @. @. \mathbb{Z}_p @. \\
@. @. @. @. @. @VVV @. \\
@. @. @. @. @. 0 @. .
\end{CD}$$
For $S\in Y(i)$ (the target) and $T\in Y(i)$ (the source), the $(S,T)$-component of $d^{''}_{i}(S,T):\Lambda^{Y(i)}\rightarrow \Lambda^{Y(i)}$ is defined as follows.

$\text{$\circ$ If $S\not=T$, $d_{i}^{''}(S,T)(x)=0$.}$

$\text{$\circ$ If $S=T$, $d_{i}^{''}(S,T)(x)=x\tau_{S}$.}$

Since $C_{\Lambda}$ is the simple complex associated to the mapping cone of  
$$d^{''}:\varprojlim_{n}\mathbb{Z}_{p}[\Gamma_{K'}/(\Gamma_{K'})^{p^{n}}]\otimes_{\mathbb{Z}_{p}}K.(\omega_{1},\ldots,\omega_{n})\rightarrow \varprojlim_{n}\mathbb{Z}_{p}[\Gamma_{K'}/(\Gamma_{K'})^{p^{n}}]\otimes_{\mathbb{Z}_{p}}K.(\omega_{1},\ldots,\omega_{n}),$$
it is quasi-isomorphic to the complex $\mathbb{Z}_{p}[[\tau]]\rightarrow \mathbb{Z}_{p}[[\tau]] :x\mapsto x\tau$, and hence to $\mathbb{Z}_{p}$.
Thus, we get the exact sequence $$0\longrightarrow\Lambda ^{X(n+1)}\mathop{\longrightarrow}^{d_n}\Lambda ^{X(n)}\mathop{\longrightarrow}^{d_{n-1}}\cdots\mathop{\longrightarrow}^{d_{0}}\Lambda ^{X(0)}\mathop{\longrightarrow}^{Aug} \mathbb{Z}_p\longrightarrow 0.$$
\end{proof}
\section{Proof of the main theorem}
\subsection{ Connection between $C_{\Gamma_{K}}(M)$ and $C_{\Lambda}$} First, let us fix some notations.
Let $G$ denote a profinite group and put $\Lambda_{G} = \mathbb{Z}_p[[G]]$. Then, $\Lambda_G$-Mod (resp.
$\mathbb{Z}_p$-Mod, $\mathscr{C}_{G}$, $\mathscr{D}_G$) denotes the category of left $\Lambda_{G}$-modules (resp. $\mathbb{Z}_p$-modules,
compact left $\Lambda_G$-modules, discrete left $\Lambda_G$-modules). Furthermore, let $D^{+}(*)$
denote the derived category of $*$ ($\in \verb+{+\Lambda_G\text{-Mod},\mathbb{Z}_p\text{-Mod}, \mathscr{C}_G,\mathscr{D}_G\verb+}+$) which consists
of complexes bounded below.

Let $M$ be a left $\Lambda$-module. Define the complex $C_{\Gamma_K}(M)$ to be
$$C_{\Gamma_K}(M) = \Hom_{\Lambda}(C_{\Lambda},M)$$
where $\Hom_{\Lambda}(A,B) \ (A, B \in \Lambda(=\Lambda_{\Gamma_{K}}$)-Mod) denotes the set of all homomorphisms $f : A \rightarrow B$
of $\Lambda$-modules. In the case $ M = D(V )$, this $C_{\Gamma_{K}}(M)$ clearly coincides with the
one defined in Section 3. On the other hand, by Proposition 4.4, we have
$$\Hom_{\Lambda}(C_{\Lambda},M)\simeq \R\Hom_{\Lambda}(\mathbb{Z}_{p},M)$$
where we denote $\R\Hom_{\Lambda}(\mathbb{Z}_{p},-) : D^{+}(\Lambda\text{-Mod}) \rightarrow  D^{+}(\mathbb{Z}_{p}\text{-Mod}).$

For every discrete left $\Lambda$-module $M$, consider the $\mathbb{Z}_{p}$-module $\Hom_{\Lambda,cont}(\mathbb{Z}_{p},M)$
of all continuous homomorphisms $f : \mathbb{Z}_{p}\rightarrow  M$ of $\Lambda$-modules. Then, we obtain
the functor
$$\Hom_{\Lambda,cont}(\mathbb{Z}_{p},-) : \mathscr{D}_{\Gamma_{K}}\rightarrow \mathscr{D}_{\mathbb{Z}_{p}}.$$
Here, $\mathscr{D}_{\mathbb{Z}_{p}}$ denotes the category $\mathscr{D}_{\verb+{+e\verb+}+}$ ($e$: unit). To define the derived functor
$\R\Hom_{\Lambda,cont}(\mathbb{Z}_{p},M)$ ($M$: discrete left $\Lambda$-module), we can use the projective resolution
of $\mathbb{Z}_{p}$ in $\mathscr{C}_{\Gamma_{K}}$ (see Remark 5.2 below). Since each component of $C_{\Lambda}$ is a  finitely
generated free $\Lambda$-module, it gives a projective resolution of $\mathbb{Z}_{p}$ in $\mathscr{C}_{\Gamma_{K}}$. Furthermore,
since we have the equality $\Hom_{\Lambda}(P,M) = \Hom_{\Lambda,cont}(P,M)$ for a finitely
generated free $\Lambda$-module $P$ and a discrete $\Lambda$-module $M$, we obtain
$$\R\Hom_{\Lambda}(\mathbb{Z}_{p},M) = \R\Hom_{\Lambda,cont}(\mathbb{Z}_{p},M).$$
If $M$ is a discrete $\Lambda$-module, we also have
$$\R\Hom_{\Lambda,cont}(\mathbb{Z}_{p},M) \simeq  \R\Gamma(\Gamma_{K},M)$$
(see [NSW, p231, (5.2.7)]). Thus, we obtain the following.
\begin{prop} If $M$ is a discrete left $\Lambda$-module, we have
$$C_{\Gamma_{K}}(M) \simeq \R\Gamma(\Gamma_{K},M).$$
\end{prop}
\begin{rem} Though it is stated in ([NSW], p231) that, to define $\R\Hom_{\Lambda_{G},cont}(L,$ $
M)$ for $L \in\mathscr{C}_{G}$ and $M \in \mathscr{D}_G$, one can use either projective resolutions of $L$ in
$\mathscr{C}_G$ or injective resolutions of $M$ in $\mathscr{D}_G$, we shall review this fact here. For a
projective resolution $P^{\cdot}\rightarrow L$ in $\mathscr{C}_G$ and an injective resolution $M \rightarrow I^{\cdot}$in $\mathscr{D}_G$, it
suffices to show
$$\Hom^{\cdot}_{\Lambda_G,cont}(L, I^{\cdot}) \rightarrow  \Hom^{\cdot}_{\Lambda_{G},cont}(P^{\cdot}, I^{\cdot})   \leftarrow \Hom^{\cdot}_{\Lambda_{G},cont}(P^{\cdot},M)$$
are quasi-isomorphisms. Here, $\Hom_{\Lambda_G,cont}(A,B)$ $(A \in \mathscr{C}_G$ and $B \in \mathscr{D}_G$) denotes
all continuous homomorphisms $f : A \rightarrow B$ of $\Lambda_{G}$-modules. For this, we have
to show that both functors $\mathscr{C}_G \rightarrow  \mathscr{D}_{\mathbb{Z}_p} : L \mapsto \Hom_{\Lambda_G,cont}(L, I)$ ($I$ is an injective
object of $\mathscr{D}_{G}$) and $\mathscr{D}_{G} \rightarrow  \mathscr{D}_{\mathbb{Z}_p} : M \mapsto \Hom_{\Lambda_G,cont}(P,M)$ ($P$ is a projective
object of $\mathscr{C}_G$) are exact functors. This follows from the fact that, for
$L \in \mathscr{C}_G$ and $M \in \mathscr{D}_G$, any continuous homomorphism $L \rightarrow M$ of $\Lambda_G$-modules
factors through a compact and discrete subgroup of $M$
.
\end{rem}
\begin{rem} The functor $C_{\Gamma_K}$ from the category $\Lambda$-Mod (resp. $\mathscr{D}_{\Gamma_K}$) to
the category $\mathbb{Z}_p$-Mod (resp. $\mathscr{D}_{\mathbb{Z}_p}$) naturally extends to the functor $C_{\Gamma_K}$ from
the derived category $D^{+}(\Lambda$-Mod) (resp. $D^{+}(\mathscr{D}_{\Gamma_K})$) to the derived category
$D^{+}(\mathbb{Z}_p$-Mod) (resp. $D^{+}(\mathscr{D}_{{Z}_p)}$). Note that the functor $C_{\Gamma_{K}}$ is an exact functor, i.e.
for an exact sequence of $\Lambda$-modules (resp. discrete $\Lambda$-modules) $0 \rightarrow  M_{1} \rightarrow M_{2} \rightarrow
M_{3} \rightarrow 0$, we have an exact sequence of complexes $0 \rightarrow C_{\Gamma_K}(M_1) \rightarrow C_{\Gamma_K}(M_2) \rightarrow C_{\Gamma_K}(M_3) \rightarrow 0$. Furthermore, Proposition 5.1 induces a canonical isomorphism
of functors $C_{\Gamma_K}(-) \simeq  \R\Gamma(\Gamma_K,-)$ from the derived category $D^{+}(\mathscr{D}_{\Gamma_K})$ to the
derived category $D^{+}(\mathscr{D}_{\mathbb{Z}_p} ).$
\end{rem}
The exact functor from the category $\mathscr{D}_{\Gamma_K}$ to the category $\Lambda$-Mod naturally
extends to the functor from the derived category $D^{+}(\mathscr{D}_{\Gamma_K})$ to the derived category
$D^{+}(\Lambda$-Mod). Therefore, the object $\R\Gamma(H_K, V )$ of the derived category
$D^{+}(\mathscr{D}_{\Gamma_K})$ gives an object of the derived category $D^{+}(\Lambda$-Mod).
\begin{prop} Let $V$ be a $p$-torsion representation of $G_K$. Then, we have an
isomorphism
$$\R\Gamma(H_K, V ) \simeq  [D(V )
\mathop{\rightarrow}^{\rho=\phi-1} D(V )]$$
in $D^{+}(\Lambda$-Mod).
\end{prop}
For the proof of this proposition, we shall introduce a subcategory of $\Lambda_{G_{K}}$-Mod
which contains the $\Lambda_{G_{K}}$-module $\mathbb{A}\otimes_{\mathbb{Z}_{p}}V$. First, let us fix some notations. Let $G$
be a profinite group and $H$ be a closed normal subgroup of $G$. Let $\mathscr{S}$ denote the
set of open subgroups of $H$ which are also normal subgroups of $G$. We define $\mathscr{E}_{G,H}$ to be the full
subcategory of $\Lambda_G$-Mod which consists of $\Lambda_G$-modules $M$ with the following
property: for all $x \in M$, there exist $U_x \in \mathscr
{S}$ and $n_x \in \mathbb{Z}_{>0}$ such that the action of
$\Ker (\Lambda_{G} \rightarrow \Lambda_{G/U_x}/p^{n_x})$ on $x$ is $0$. Then, $\mathscr{E}_{G,H}$ forms an abelian category.
\begin{lem} The category $\mathscr
{E}_{G,H}$ has sufficiently many injectives.
\end{lem}
\begin{proof} For $M \in \mathscr{E}_{G,H}$, there exists an inclusion $M \hookrightarrow I$ where $I$ is an injective
object of $\Lambda_{G}$-Mod. Define $I'$ to be $\verb+{+x \in I\verb+|+\exists U \in\mathscr{ S}, n \in \mathbb{Z}_{>0}$ s.t. the action of
$\Ker (\Lambda_{G} \rightarrow  \Lambda_{G/U}/p^n)$ on $x$ is $0$ $\verb+}+$. Then, $I'$ becomes an injective object of $\mathscr{E}_{G,H}$
such that $M \subset I'$.
\end{proof}
\begin{lem} \begin{enumerate}
\item For $U, U'\in \mathscr{S}$, $U' \subset U$, the homomorphism $\Lambda_{G/U'}\otimes_{\Lambda_{H/U'}}\Lambda_{H/U}\rightarrow \Lambda_{G/U}$ is an isomorphism.

\item For $U\in \mathscr{S}$, $\Lambda_{G/U}$ is flat as a right $\Lambda_{H/U}$-module.
\end{enumerate}
\end{lem}
\begin{proof} (1) The natural homomorphism $G/U' \rightarrow G/H$ has a continuous section
$s : G/H \rightarrow G/U'$ (see [S2], p4, Proposition 1.). With this, we obtain a homeomorphism
$G/H \times H/U' \simeq G/U' : (a, b)\mapsto  s(a) \cdot b$ of profinite sets which is
compatible with the right action of $H/U'$. Therefore, we get an isomorphism
$f' : \mathbb{Z}_{p}[[G/H]]\otimes_{\mathbb{Z}_{p}}\mathbb{Z}_{p}[H/U'] \simeq \mathbb{Z}_{p}[[G/U']]$ of right $\mathbb{Z}_{p}[H/U']$-modules. By using
the composition with the section $s$ and $G/U' \rightarrow G/U$, we similarly get an isomorphism
$f : \mathbb{Z}_p[[G/H]]\otimes_{\mathbb{Z}_{p}}\mathbb{ Z}_p[H/U] \simeq \mathbb{Z}_{p}[[G/U]]$ of right $\mathbb{Z}_{p}[H/U]$-modules. Since
$f$ and $f'$ are compatible with $\mathbb{Z}_p[H/U'] \rightarrow \mathbb{Z}_{p}[H/U]$ and $\mathbb{Z}_p[[G/U']] \rightarrow \mathbb{Z}_{p}[[G/U]]$,
we obtain the desired result.

(2) Since we have the isomorphism $f : \mathbb{Z}_p[[G/H]]\otimes_{\mathbb{Z}_{p}}\mathbb{Z}_{p}[H/U]\simeq \mathbb{Z}_{p}[[G/U]]$ of
right $\mathbb{Z}_p[H/U]$-modules and $\mathbb{Z}_p[[G/H]]$ is flat as a $\mathbb{Z}_p$-module, $\Lambda_{G/U}$ is flat as a
right $\Lambda_{H/U}$-module.
\end{proof}
For $M \in \mathscr{D}_H$ and $U \in \mathscr{S}$, define $M^U=\verb+{+x \in M\verb+|+\text{the action of} \ \Ker (\Lambda_{H} \rightarrow \Lambda_{H/U} ) \  
\text{on $x$ is trivial}\verb+}+$. Since $M$ is an object of $\mathscr{D}_{H}$, we have $M = \varinjlim _{U\in\mathscr{S}}M^U$. Define
the left $\Lambda_{G/U}$-module $T_{U}(M)$ to be $\Lambda_{G/U}\otimes_{\Lambda_{H/U}}M^{U}$. 
By Lemma 5.6.(1),
for $U' \in\mathscr{ S}$, $U' \subset U$, the natural morphism $\Lambda_{G/U'}\otimes_{\Lambda_{H/U'}}M^{U} \rightarrow T_{U}(M)$ becomes
an isomorphism. Therefore, by Lemma 5.6.(2), we obtain an injection
$T_{U}(M)\rightarrow T_{U'}(M)$ which is compatible with the action of $\Lambda_G$. Then, it follows
easily that $\verb+{+\Lambda_{U}(M)\verb+|+U \in\mathscr{ S}\verb+}+$ forms an inductive system. We denote the inductive
limit $\varinjlim _{U\in\mathscr{S}}T_{U}(M)$ by $T(M)$. Since $T(M)$ becomes an object of $\mathscr{E}_{G,H}$, we obtain
a functor $T : \mathscr{D}_H \rightarrow \mathscr{E}_{G,H}$. Furthermore, by Lemma 5.6.(2) and the fact
$M = \varinjlim _{U\in\mathscr{S}}M^U$, it follows that the functor $T$ is an exact functor.
\begin{lem} If $H$ is a finite group, $\Ker (\Lambda_{G} \rightarrow \Lambda_{G/H})$ is generated by $\verb+{+h-1\verb+|+h \in H\verb+}+$.
\end{lem}
\begin{proof} There exists an exact sequence of projective systems of finite abelian
groups
$$\bigoplus_{h\in H} \mathbb{Z}/p^n[G/V ] \cdot (h - 1) \rightarrow \mathbb{Z}/p^n[G/V ] \rightarrow \mathbb{Z}/p^n[G/(V \cdot H)] \rightarrow 0$$
where $n$ and $V$ run through positive integers and open normal subgroups of $G$.
Since these are projective systems of finite abelian groups, the filtered projective
limit preserves the exactness by Pontryagin duality. Thus, we obtain an exact sequence
$$\bigoplus_{h\in H}
\Lambda_{G} \cdot (h - 1) \rightarrow \Lambda_{G} \rightarrow \Lambda_{G/H} \rightarrow 0.$$
\end{proof}
\begin{lem} Let $N$ be an object of $\mathscr{E}_{G,H}$. If the action of $\Ker (\Lambda_{H} \rightarrow \Lambda_{H/U} )$ on
$x$ is $0$ for $x \in N, U \in \mathscr{S}$, then, the action of $\Ker (\Lambda_G \rightarrow \Lambda_{G/U})$ on $x$ is also $0$.
\end{lem}
\begin{proof} By the definition of $\mathscr{E}_{G,H}$, there exists an element $U' \in\mathscr{S}$ contained in $U$
such that the action of $\Ker (\Lambda_G \rightarrow \Lambda_{G/U'})$ on $x$ is $0$. By applying Lemma 5.7 above
to $U/U' \subset G/U'$, we see that $\Ker (\Lambda_{G/U'} \rightarrow \Lambda_{G/U})$ is an ideal generated by $\verb+{+g-1\verb+|+g \in U/U'\verb+}+$.
Since the action of $U$ on $x$ is trivial by hypothesis, the action of this ideal on $x$
is $0$. 
\end{proof}
\begin{prop} The functor $T$ is a left-adjoint functor of the forgetful functor
$F : \mathscr{E}_{G,H }\rightarrow \mathscr{D}_H$.
\end{prop}
\begin{proof} For an object $M$ of $\mathscr{D}_H$, the natural map $M^U \rightarrow T_U(M) : x \mapsto 1\otimes 
 x$ is a
homomorphism of $\Lambda_H$-modules and compatible with respect to $U$. By taking the
inductive limit, we obtain $\alpha_M : M \rightarrow F \circ T(M)$. This morphism is functorial in
$M$. On the other hand, for an object $N$ of $\mathscr{E}_{G,H}$, $N^U$ becomes a $\Lambda_{G/U}$-module by
Lemma 5.8 above. Therefore, we have a homomorphism $T_U(N) \rightarrow N^U$ of $\Lambda_{G/U}$-modules
and this homomorphism is compatible with respect to $U$. By taking the
inductive limit, we obtain $\beta_{N} : T \circ F(N) \rightarrow N$. This morphism is functorial in 
$N$. For $M \in \mathscr{D}_{H}$ and $N \in \mathscr{E}_{G,H}$, we obtain maps which are functorial in $M$ and
$N$
\begin{align*}\Hom_{\mathscr{E}_{G,H}} (T(M),N) \rightarrow \Hom_{\mathscr{D}_H}(M, F(N)) : \varphi \mapsto F(\varphi) \circ \alpha_{M},\\
\Hom_{\mathscr{D}_{H}}(M, F(N)) \rightarrow \Hom_{\mathscr{E}_{G,H}} (T(M),N) : \psi  \mapsto \beta_{N}\circ T(\psi ).
\end{align*}
It follows easily that each map is inverse to the other map. 
\end{proof}
Since the functor $T$ is exact and a left-adjoint functor of $F$ by Proposition 5.9,
the functor $F$ preserves injective objects.

Now, for an object $N$ of $\mathscr{E}_{G,H}$, define $N^H = \verb+{+x \in N\verb+|+h(x) = x, \forall h \in H\verb+}+$.
\begin{lem} $N^H$ is a left $\Lambda_{G/H}$-module
\end{lem}
\begin{proof} For $x \in N^H$, there exists an element $U \in \mathscr{S}$ such that the action of
$\Ker (\Lambda_{G} \rightarrow \Lambda_{G/U})$ on $x$ is $0$. By applying Lemma 5.7 to $H/U \subset G/U$, it follows
that the ideal $\Ker (\Lambda_{G/U} \rightarrow \Lambda_{G/H})$ is generated by $\verb+{+h -1\verb+|+ h \in H/U\verb+}+$. Thus, we
see that the action of this kernel on $x$ is $0$. 
\end{proof}
With this, we have a left exact functor $\Gamma_{\mathscr{E}}(H,-) : \mathscr{E}_{G,H} \rightarrow \Lambda_{G/H}$-Mod $: N \mapsto N^{H}$ and
$$\R\Gamma_{\mathscr{E}}(H,-) : D^{+}(\mathscr{E}_{G,H}) \rightarrow D^{+}(\Lambda_{G/H}\text{-Mod}).$$
\begin{prop} The following diagram is commutative
$$\begin{CD}
 D^{+}(\mathscr{E}_{G,H}) @>\R\Gamma_{\mathscr{E}}(H,-)>> D^{+}(\Lambda_{G/H}\text{-Mod}) \\
@VF_{1}VV   @VF_{2}VV \\
D^{+}(\mathscr{D}_{H}) @>\R\Gamma(H,-)>>D^{+}(\mathbb{Z}_{p}\text{-Mod}).   
\end{CD}$$ 
Here the two vertical arrows denote the functors induced by the forgetful functors $\mathscr{E}_{G,H}\rightarrow  \mathscr{D}_{H}$
and $\Lambda_{G/H}\text{-Mod}\rightarrow  \mathbb{Z}_{p}\text{-Mod}.$
\end{prop}
\begin{proof} We have a commutative diagram
$$\begin{CD}
 \mathscr{E}_{G,H} @>\Gamma_{\mathscr{E}}(H,-)>>\Lambda_{G/H}\text{-Mod} \\
@VVV   @VVV \\
\mathscr{D}_{H} @>\Gamma(H,-)>>   \mathbb{Z}_{p}\text{-Mod}.
\end{CD}$$ 
The two vertical functors are exact and the left vertical map preserves injective objects
by Proposition 5.9. Thus, it follows easily that the diagram in this proposition
is commutative. 
\end{proof}
\begin{prop} Let $F_3$ (resp. $F_4$) be the functor $D^{+}(\mathscr{D}_{G}) \rightarrow  D^{+}(\mathscr{E}_{G,H})$ (resp.
$D^{+}(\mathscr{D}_{G/H}) \rightarrow D^{+}(\Lambda_{G/H}\text{-Mod}))$ induced by the inclusion functor $\mathscr{D}_G \rightarrow \mathscr{ E}_{G,H}$
(resp. $\mathscr{D}_{G/H}\rightarrow \Lambda_{G/H}\text{-Mod})$. Then, the following diagram is commutative
$$\begin{CD}
D^{+}(\mathscr{D}_{G})  @>\R\Gamma(H,-)>>  D^{+}(\mathscr{D}_{G/H})\\
@VF_{3}VV   @VF_{4}VV \\
D^{+}(\mathscr{E}_{G/H})@>\R\Gamma_{\mathscr{E}}(H,-)>> D^{+}(\Lambda_{G/H}\text{-Mod}).       
\end{CD}$$ 
\end{prop}
\begin{proof}It suffices to show that, for an injective object $I$ of $\mathscr{D}_{G}$, we have $\R^{i}\Gamma_{\mathscr{E}}(H,$ $
F_3(I)) = 0 \ (i > 0)$. By Proposition 5.11, we have an isomorphism $\R^{i}\Gamma_{\mathscr{E}}(H, F_3(I))$
$= \R^{i}\Gamma(H, F_1 \circ F_3(I))$ of $\mathbb{Z}_p$-modules. Since the following diagram is commutative
by the group cohomology theory for discrete modules, we obtain $\R^{i}\Gamma(H, F_1 \circ F_3(I)) = F_2 \circ F_4( \R^i\Gamma(H, I)) = 0.$
$$\begin{CD}
  D^{+}(\mathscr{D}_{G}) @>\R\Gamma(H,-)>> D^{+}(\mathscr{D}_{G/H})\\
@VF_{1}\circ F_{3}VV   @VF_{2}\circ F_{4}VV \\
D^{+}(\mathscr{D}_{H})@>\R\Gamma(H,-)>> D^{+}(\mathbb{Z}_{p}\text{-Mod}).  
\end{CD}$$ 
\end{proof}
Now, we shall give the proof of Proposition 5.4. Note that, since $\mathbb{A}\otimes_{\mathbb{Z}_{p}}V$ 
becomes an object of $\mathscr{E}_{G_K,H_K}$ (see Remark 2.6), we have an exact sequence
$$0\rightarrow V\rightarrow \mathbb{A}\otimes_{\mathbb{Z}_{p}}V\mathop{\rightarrow}^{\rho=\phi-1} \mathbb{A}\otimes_{\mathbb{Z}_{p}}V\rightarrow 0$$
in $\mathscr{E}_{G_K,H_K}$. First, we will show that we have
$$H^i(H_K,\mathbb{A}\otimes_{\mathbb{Z}_{p}}V) = 0 \quad \text{for all $i > 0$}.$$
Since we have the canonical isomorphism of Galois groups $H_{K} \simeq G_{\mathbb{E}_{K}}$ by the theory
of field of norms, we have only to show $H^{i}(G_{\mathbb{E}_{K}},\mathbb{A}\otimes_{\mathbb{Z}_{p}}V)=0$ 
for all $i > 0$.
On the other hand, we have isomorphisms of $G_{\mathbb{E}_{K}} \ (\simeq  H_K)$-modules $\mathbb{A}\otimes_{\mathbb{Z}_{p}}V\simeq \mathbb{A}\otimes_{\mathbb{A}_{K}}D(V)\simeq\bigoplus_{j=1}^{d}\mathbb{A}/p^{m_{j}}\mathbb{A}$. 
Thus, it suffices to show $H^i(G_{\mathbb{E}_K},\mathbb{A}/p^m\mathbb{A}) =
0$ for all $i > 0.$ This is clear for $m = 1$ ($H^i(G_{\mathbb{E}_K},\mathbb{E}) = 0$ for all $i > 0$) and
the general case can be deduced by induction on the integer $m$. Thus, by using
Proposition 5.11, we obtain isomorphisms in $D^{+}(\Lambda\text{-Mod})$ from the exact
sequence above
\begin{align*}\R\Gamma_{\mathscr{E}}(H_K, V ) \simeq &\R\Gamma_{\mathscr{E}}(H_K, [\mathbb{A}\otimes_{\mathbb{Z}_{p}}V\mathop{\rightarrow}^{\phi-1} \mathbb{A}\otimes_{\mathbb{Z}_{p}}V])\\
\simeq & \Gamma_{\mathscr{E}}(H_{K},[\mathbb{A}\otimes_{\mathbb{Z}_{p}}V\mathop{\rightarrow}^{\phi-1} \mathbb{A}\otimes_{\mathbb{Z}_{p}}V])\\
=&[D(V)\mathop{\rightarrow}^{\phi-1} D(V)].
\end{align*}
On the other hand, by Proposition 5.12, $\R\Gamma_{\mathscr{E}}(H_K, V )$ coincides with the image of
the Galois cohomology $\R\Gamma(H_K, V ) \in D^{+}(\mathscr{D}_{\Gamma_{K}})$ by the functor $F_4 : D^{+}(\mathscr{D}_{\Gamma_{K}}) \rightarrow D^{+}(\Lambda\text{-Mod})$. Thus, this completes the proof of Proposition 5.4.
\subsection{Conclusion} We now compute the Galois cohomology $\R\Gamma(G_K, V )$ for a $p$-torsion
representation of $V$ of $G_K$. We have
$$\R\Gamma(G_K, V ) \simeq \R\Gamma(\Gamma_K, \R\Gamma(H_K, V )).$$
From Proposition 5.1 and Remark 5.3, we obtain
$$\R\Gamma(\Gamma_K, \R\Gamma(H_K, V )) \simeq C_{\Gamma_K}(\R\Gamma(H_K, V )).$$
By Proposition 5.4,
$$\R\Gamma(G_K, V ) \simeq C_{\Gamma_K}([D(V )\mathop{\rightarrow}^{\rho}D(V )]) \simeq C_{\phi, \Gamma_{K}}(D(V )).$$
Thus, this completes the proof of the main theorem.

{\bf Acknowledgment} The author is grateful to his advisor Professor Kazuya Kato for his continuous advice and encouragements.  He is also grateful to the referee for careful reading and numerous detailed and helpful comments. 
A part of this work was done while he was staying at Universit\'e Paris-Sud 11 and
he thanks this institute for the hospitality.
His staying at Universit\'e Paris-Sud 11 is partially
supported by JSPS Core-to-Core Program
``New Developments of Arithmetic Geometry, Motive, Galois Theory, and Their Practical Applications''
and he thanks Professor Makoto Matsumoto for encouraging this visiting.
This research is partially
supported by JSPS Research Fellowships for Young Scientists.

{\bf References}
[A]   Andreatta, F.: Generalized ring of norms and generalized $(\phi,\Gamma)$-modules. Preprint.

[Bn]  Benois, D.: On Iwasawa theory of crystalline representations. Duke. Math. 104 (2000), 211--267.

[Br1]  Berger, L.: Repr\'esentations $p$-adiques et \'equations diff\'erentielles. Invent. Math. 148 (2002), 219--284. 

[Br2]  Berger, L.: An introduction to the theory of $p$-adic representations. Geometric aspects of Dwork theory. Vol. I, II, 255--292, 2004. 

[F]   Fontaine, J-M.: Repr\'esentations $p$-adiques des corps locaux. I. The Grothendieck Festschrift, Vol. II, 249--309, Progr. Math., 87, Birkhauser, 1990. 

[FW]  Fontaine, J-M., Wintenberger, J-P.: Le "corps des normes" de certaines extensions alg\'ebriques de corps locaux. C. R. Acad. Sci. Paris Ser. A-B 288 (1979), A367--A370. 

[H1]  Herr, L.: Sur la cohomologie galoisienne des corps $p$-adiques. Bull. Soc. Math. France 126 (1998), 563--600.

[H2]  Herr, L.: $\Phi$-$\Gamma$-modules and Galois cohomology. Invitation to higher local fields (Munster, 1999), 263--272 (electronic), Geom. Topol. Monogr., 3, Geom. Topol. Publ., Coventry, 2000. 

[H3]  Herr, L.: Une approche nouvelle de la dualit\'e locale de Tate. Math. Ann. 320 (2001), 307--337. 

[NSW]  Neukirch, J., Schmidt, A., Wingberg, K.: Cohomology of number fields. Grundlehren der Mathematischen Wissenschaften [Fundamental Principles of Mathematical Sciences], 323. Springer-Verlag, Berlin, 2000. 

[S1]  Serre, J-P.: Local fields. Graduate Texts in Mathematics, 67. Springer-Verlag, New York-Berlin, 1979. 

[S2] Serre, J-P.: Galois cohomology. Springer-Verlag, Berlin, 1997 

[W]  Wintenberger, J-P.: Le corps des normes de certaines extensions infinies de corps locaux; applications. Ann. Sci. \'Ecole Norm. Sup. (4) 16 (1983), 59--89.

\end{document}